\documentclass[12pt]{amsart}

\usepackage[utf8]{inputenc}
\usepackage{amssymb, latexsym, amsmath, amsfonts}

\newtheorem{theorem}{Theorem}
\newtheorem{corollary}[theorem]{Corollary}

\usepackage{bbm}
\usepackage{xcolor}
\usepackage{hyperref}

\title[Mean values of the Riemann zeta function]{Mean values of the Riemann zeta function at shifted zeros under the Riemann Hypothesis}

\begin{document}

\author[R. Garunk\v{s}tis]{Ram\={u}nas Garunk\v{s}tis}
\address[Ram\={u}nas Garunkštis]{Institute of Mathematics\\Faculty of Mathematics and Informatics\\Vilnius University\\Naugarduko 24, LT-03225, Vilnius\\Lithuania}
\email{ramunas.garunkstis@mif.vu.lt}

\author[J. Paliulionyt\.{e}]{Julija Paliulionyt\.{e}}
\address[Julija Paliulionyt\.{e}]{Institute of Mathematics\\Faculty of Mathematics and Informatics\\Vilnius University\\Naugarduko 24, LT-03225, Vilnius\\Lithuania}
\email{julija.paliulionyte@mif.stud.vu.lt}

\begin{abstract}
Assuming the Riemann hypothesis, we obtain asymptotic formulas for $\sum_{0<\gamma<T}\zeta(\rho+\delta)\zeta(1-\rho+\overline{\delta})$ in the region $-\frac{a}{\log T} \leq \Re \delta \leq \frac{1}{2}+\frac{a}{\log T}$, $|\Im \delta|\ll 1$. Unconditionally, this asymptotic formula was recently obtained by Garunk\v{s}tis and Novikas in essentially the same region, with a slight incompleteness. Assuming RH, we obtain a sharper error term, and we also correct an inaccuracy in the unconditional error term there.
\end{abstract}

\maketitle

\section{Introduction}

Let $\rho=\beta+i\gamma$ denote a nontrivial zero of $\zeta(s)$. For $\delta=\delta_1+i\delta_2 \in \mathbb{C}$ define
\begin{equation*}
    \mathcal{S}(\delta)=\mathcal{S}(\delta,T)=\sum_{0<\gamma<T} \zeta(\rho+\delta)\zeta(1-\rho+\bar{\delta}).
\end{equation*}

For real $\alpha$,   $|\alpha|\le\frac{1}{4\pi}\log\frac{T}{2\pi}$, assuming the Riemann hypothesis (RH), Gonek \cite{Gonek} proved
\begin{equation*}
{\mathcal S}\left(\frac{2\pi i\alpha} {\log \frac{T}{2\pi}}\right) = \left(1-\left(\frac{\sin(\pi\alpha)}{\pi \alpha}\right)^2\right)\frac{T}{2\pi}\log^2 T+O(T\log T)
\end{equation*}
uniformly in $\alpha$.

Garunk\v{s}tis and Novikas \cite{Garunkstis} considered unconditional asymptotic formulas for $\mathcal{S}(\delta)$ which are uniform in $-\frac{a}{\log T} \leq \delta_1 \leq \frac{1}{2}+\frac{a}{\log T}$, $|\delta_2|\ll 1$, where $a>0$. We refer to their paper for a historical overview of related results. To that overview, we add the following two references. Heap, Li, and Zhao \cite[Theorem 2]{hlz22} obtained
\begin{align}\label{hlz}
\sum_{0\le \gamma\le T}\zeta(\rho+\alpha)\zeta(1-\rho+\beta)Q(\rho)\overline{Q}(1-\rho)
  =  main\ term+O(T(\log T)^{-A}),
\end{align}
where $\alpha,\beta \ll 1/\log T$, $\alpha,\beta \in \mathbb C$, $Q(s)$ is a certain Dirichlet polynomial, and $A$ is an arbitrary positive constant. For $Q\equiv 1$, $\alpha=\delta$, and $\beta=\overline{\delta}$ their main term coincides with the asymptotics in \eqref{S(delta)} below up to the error $O(T^{1/2+\epsilon})$. The authors note that, assuming the Generalized Riemann hypothesis, the error term in \eqref{hlz} can be replaced by $O(T^{1/2+\epsilon})$.  
Recently Benli, Elma, and Ng \cite[Theorem 1.1]{ben24} studied 
\begin{align*}
          \sum_{0\le \gamma\le T}\zeta(\rho+\alpha)Q(\rho)P(1-\rho)
\end{align*}
where $P$ is also a Dirichlet polynomial.

In \cite{Garunkstis} it is shown that if $\delta_1\neq 0$ and 
\begin{align}\label{0infty}
|\delta|\sqrt{T}\log^2T \to 0 \quad \text{or}\quad|\delta|T^{\frac{A}{2}(\log T)^{-\frac{2}{5}}(\log\log T)^{-\frac{1}{5}}} \to \infty, 
\end{align} 
then
\begin{align}\label{S(delta)}
\mathcal{S}(\delta) 
\sim&
\zeta(1+2\delta_1)\frac{T}{2\pi}\log{\frac{T}{2\pi e}}+2\Re\frac{\zeta'}{\zeta}(1+\delta)\zeta(1+2\delta_1)\frac{T}{2\pi}\nonumber
\\& -
2\Re\zeta(1+\overline{\delta})\frac{\zeta(1-\delta)}{1-\delta}\left(\frac{T}{2\pi}\right)^{1-\delta}
\\&+
2\chi(\delta_1)\Re\frac{\zeta'}{\zeta}(1-\overline{\delta})\frac{\zeta(1-2\delta_1)}{1-2\delta_1}\left(\frac{T}{2\pi}\right)^{1-2\delta_1}\nonumber
\\&+
\chi(\delta_1)\frac{\zeta(1-2\delta_1)}{1-2\delta_1}\left(\frac{T}{2\pi}\right)^{1-2\delta_1}\left(\log{\frac{T}{2\pi}}-\frac{1}{1-2\delta_1}\right),\nonumber
\end{align}
where
\begin{equation*}
    \chi(\delta_1)=
    \begin{cases}
        0 \quad \text{if} \quad \delta_1 > (\log T)^{-\frac25}(\log\log T)^{-\frac15},\\
        1 \quad \text{if} \quad \delta_1 \leq (\log T)^{-\frac25}(\log\log T)^{-\frac15}.
    \end{cases}
\end{equation*}
Taking a limit $\delta_1 \to 0$ we get an asymptotic for $\mathcal{S}(\delta)$ when $\delta_1=0$, $\delta_2\ \neq 0$. Note that the behavior of $\mathcal{S}(\delta)$ changes depending on whether $T^{-\delta_1}$ tends to zero, stays constant, or tends to 1  as $T$ increases (see \cite[Section 2]{Garunkstis}). By \eqref{0infty} we see that asymptotic \eqref{S(delta)} wasn't proved in entire region $-\frac{a}{\log T} \leq \delta_1 \leq \frac{1}{2}+\frac{a}{\log T}$, $|\delta_2|\ll 1$, for example it wasn't proved in $T^{-\frac12}\log^{-2}T\le|\delta|\le T^{-A\log^{-\frac25} T(\log\log T)^{-\frac15}/2}$.

In \cite{Garunkstis}, the error term $E(T)$ in the asymptotic \eqref{S(delta)} requires a correction,  which is given in Section~\ref{correction}.  The correction doesn't affect the above statements concerning $\mathcal{S}(\delta)$.

We have that for all $\epsilon>0$ the asymptotic \eqref{S(delta)} holds in $|\delta|\log T \geq \epsilon$, $\delta_1 \neq 0$. In this paper, we will show that, assuming the Riemann hypothesis, the asymptotic \eqref{S(delta)} holds in $|\delta|\log T \leq \epsilon_0$, $\delta_1 \neq 0$ for some $\epsilon_0>0$, and hence holds everywhere in the region $-\frac{a}{\log T} \leq \delta_1 \leq \frac{1}{2}+\frac{a}{\log T}$, $|\delta_2|\ll 1$, $\delta_1 \neq 0$.

In Section \ref{sec3}, following the proof of Theorem 3 in \cite{Garunkstis}, we obtain
\begin{theorem}\label{thm1} Assume the RH and that $\delta \log T \to 0$. Then, for every $\epsilon > 0$,
\begin{equation*}
\mathcal{S}(\delta)=\frac{|\delta|^2}{12}\frac{T}{2\pi}\log^4T+O\left( |\delta|^3T^{1+\epsilon}+|\delta|^2T\log^3T\right).
\end{equation*}
\end{theorem}
This theorem gives an asymptotic if $|\delta|\leq T^{-\epsilon}$. If $\delta \log T \to 0$, then the Taylor series expansion shows that the right-hand side of $\eqref{S(delta)}$ is also asymptotic to $\frac{|\delta|^2}{12}\frac{T}{2\pi}\log^4T$ (see \cite[Section 6]{Garunkstis}).

In Section \ref{sec4}, following the proof of Theorem 1 in \cite{Garunkstis}, we show
\begin{theorem}\label{thm2} Assume the RH. Let $-\frac{a}{\log T} \leq \delta_1 \leq \frac12+\frac{a}{\log T}$, $\delta_1 \neq 0$, and $|\delta_2|\ll 1$. Then, for every $\epsilon>0$,
\begin{align}\label{Sthm2}
    \mathcal{S}(\delta)=&\zeta(1+2\delta_1)\frac{T}{2\pi}\left(\log\frac{T}{2\pi}-1\right)\nonumber
    \\&+
    \mathbbm{1}_{\left[-\frac{a}{\log T},\frac14\right)}(\delta_1)\frac{\zeta(1-2\delta_1)}{1-2\delta_1}\left(\frac{T}{2\pi}\right)^{1-2\delta_1}\left(\log\frac{T}{2\pi}-\frac{1}{1-2\delta_1}\right)\nonumber
    \\&+
    2\Re\left(\frac{\zeta'}{\zeta}(1+\delta)\zeta(1+2\delta_1)\frac{T}{2\pi}\right)
    \\&+
    2\mathbbm{1}_{\left[-\frac{a}{\log T},\frac14\right)}(\delta_1)\Re\left(\frac{\zeta'}{\zeta}(1-\overline{\delta})\frac{\zeta(1-2\delta_1)}{1-2\delta_1}\left(\frac{T}{2\pi}\right)^{1-2\delta_1}\right)\nonumber
    \\&
    -2\Re\left(\zeta(1+\overline{\delta})\frac{\zeta(1-\delta)}{1-\delta}\left(\frac{T}{2\pi}\right)^{1-\delta}\right)\nonumber
    +
    O\left(T^{\frac{1}{2}+\epsilon}\right).\nonumber
\end{align}
\end{theorem}

Finally, in Section \ref{sec3} we prove that there is $\epsilon_0>0$ such that the main term of the formula \eqref{Sthm2} dominates uniformly in $|\delta| \log T \leq \epsilon_0$, $|\delta| \geq T^{-\frac14+\epsilon}$, $\delta_1 \neq 0$ which leads to the following corollary.
\begin{corollary}\label{cor}
Assume the RH. Let $-\frac{a}{\log T} \leq \delta_1 \leq \frac{1}{2}+\frac{a}{\log T}$ and $|\delta_2|\ll 1$.
Then for $\mathcal{S}(\delta)$ the asymptotic \eqref{S(delta)} is true. The case $\delta_1=0$ in \eqref{S(delta)} is understood as a limit.
\end{corollary}

\section{Error term $E(T)$}\label{correction}

In this section, we don't assume RH. The error term in \cite{Garunkstis} was
\begin{align*}
    E(T)=
    \begin{cases}
        T^{1-A(\log T)^{-\frac25} (\log\log T)^{-\frac15}} &\text{if} \quad -\frac{a}{\log T} \leq \delta_1 \leq (\log T)^{-\frac25} (\log\log T)^{-\frac15},\\
        T^{1-A(\delta_1\log T)^{-\frac23}(\log\log T)^{-\frac13}-\delta_1} &\text{if} \quad (\log T)^{-\frac25} (\log\log T)^{-\frac15} \leq \delta_1 \leq \frac12+\frac{a}{\log T},
    \end{cases}
\end{align*}
where constant $A=1/(4\cdot 57.54)$ is related to the zero free region of $\zeta(s)$ (\cite[formula (5)]{Garunkstis}). However, the second case for $E(T)$ is not completely correct. The mistake in the proof was in the bound (see \cite[page 214]{Garunkstis}) that stated that
\begin{align*}
&\int_{d-iU}^{d+iU}\frac{\zeta'}{\zeta}(s+\delta)\zeta(s+2\delta_1)\zeta(s)\left(\frac{T}{2\pi}\right)^s\frac{d s}{s}\\
&\ll
T^{1-2A(\log U)^{-\frac23}(\log\log U)^{-\frac13}-\delta_1}\log^4 T,
\end{align*}
where  $U\ll T^{1-\varepsilon}$ and 
\begin{align*}
d= 1-3A\log^{-\frac23} U(\log\log U)^{-\frac13}-\delta_1,
\end{align*}
if 
\begin{align*}
|\delta_1-3A\log^{-\frac23} U(\log\log U)^{-\frac13}|>A\log^{-\frac23} U(\log\log U)^{-\frac13},
\end{align*}
 otherwise  
\begin{align*}
d= 1-2A\log^{-\frac23} U(\log\log U)^{-\frac13}-\delta_1.
\end{align*}
We will show that
\begin{align*}
    E(T)=
    \begin{cases}
        T^{1-A(\log T)^{-\frac25} (\log\log T)^{-\frac15}} &\text{if} \quad -\frac{a}{\log T} \leq \delta_1 \leq (\log T)^{-\frac25} (\log\log T)^{-\frac15},\\
        T^{1-A(\delta_1\log T)^{-\frac23}(\log\log T)^{-\frac13}-\delta_1+\delta_1^2/2} &\text{if} \quad (\log T)^{-\frac25} (\log\log T)^{-\frac15} \leq \delta_1 \leq \frac12+\frac{a}{\log T}.
    \end{cases}
\end{align*}
For $2 \leq |t|\leq U$, we have (see \cite[formulas (22) and (43)]{Garunkstis})
\begin{align*}
&\frac{\zeta'}{\zeta}(d+it+\delta)\ll \log T,\\
&\zeta(d+it+2\delta_1)\ll \log T,\\
&\zeta(d+it)\ll t^{\frac12(1-d)}\log t\leq T^{\frac32A(\log U)^{-\frac23}(\log\log U)^{-\frac13}+\frac{\delta_1}{2}}\log T.
\end{align*}
Hence,
\begin{align}\label{integral}
&\int_{d-iU}^{d+iU}\frac{\zeta'}{\zeta}(s+\delta)\zeta(s+2\delta_1)\zeta(s)\left(\frac{T}{2\pi}\right)^s\frac{d s}{s}\\
&\ll
U^{\frac32A(\log U)^{-\frac23}(\log\log U)^{-\frac13}+\frac{\delta_1}{2}}
T^{1-2A(\log U)^{-\frac23}(\log\log U)^{-\frac13}-\delta_1}\log^3 T.\nonumber
\end{align}
Take
\begin{equation*}
U=T^{(\log T)^{-\frac25}(\log\log T)^{-\frac15}+\delta_1-\delta_1^2/2}.
\end{equation*}
We will prove that the right-hand side of \eqref{integral} is $\ll E(T)$.

\textbf{Case 1}: $-\frac{a}{\log T} \leq \delta_1 \leq (\log T)^{-\frac25} (\log\log T)^{-\frac15}$.

\noindent We have
\begin{align*}
    &U \leq T^{(\log T)^{-\frac25}(\log\log T)^{-\frac15}+\delta_1} \leq T^{2(\log T)^{-\frac25}(\log\log T)^{-\frac15}}\\
    &\implies \log U \leq 2(\log T)^{\frac35}(\log\log T)^{-\frac15}\\
    &\implies (\log U)^{-\frac23} \geq 2^{-\frac23}(\log T)^{-\frac25}(\log\log T)^{\frac{2}{15}}.
\end{align*}
Also,
\begin{equation*}
    U \leq T \implies \log \log U \leq \log \log T \implies (\log \log U)^{-\frac13} \geq (\log \log T)^{-\frac13}.
\end{equation*}
Hence,
\begin{align}\label{T}
    T^{1-2A(\log U)^{-\frac23}(\log\log U)^{-\frac13}-\delta_1} &\leq T^{1-2^{\frac13}A(\log T)^{-\frac25}(\log \log T)^{-\frac15}-\delta_1} \nonumber
    \\&\ll
    T^{1-2^{\frac13}A(\log T)^{-\frac25}(\log \log T)^{-\frac15}}.
\end{align}
Moreover,
\begin{align}\label{U}
    U^{\frac32A(\log U)^{-\frac23}(\log\log U)^{-\frac13}+\frac{\delta_1}{2}} &\leq e^{\frac32A(\log U)^{\frac13}}U^{\frac12 (\log T)^{-\frac25}(\log \log T)^{-\frac15}} \nonumber
    \\&\leq
    e^{3A2^{-\frac23}(\log T)^{\frac15}(\log\log T)^{-\frac{1}{15}}}T^{(\log T)^{-\frac45}(\log \log T)^{-\frac25}} \nonumber
    \\&\leq
    T^{3A(\log T)^{-\frac45}(\log\log T)^{-\frac{1}{15}}+(\log T)^{-\frac45}(\log \log T)^{-\frac25}}
    \\&=
    T^{o\left((\log T)^{-\frac25}(\log \log T)^{-\frac15}\right)}. \nonumber
\end{align}
By \eqref{T} and \eqref{U} we find that the right-hand side of \eqref{integral} is $\ll T^{1-A(\log T)^{-\frac25}(\log \log T)^{-\frac15}}$.

\textbf{Case 2}: $(\log T)^{-\frac25}(\log\log T)^{-\frac15} \leq \delta_1 \leq \frac12+\frac{a}{\log T}$.

\noindent We have
\begin{align*}
    &U \leq T^{(\log T)^{-\frac25}(\log\log T)^{-\frac15}+\delta_1} \leq T^{2\delta_1} \implies \log U \leq 2\delta_1\log T\\
    &\implies (\log U)^{-\frac23} \geq 2^{-\frac23}(\delta_1\log T)^{-\frac23}.
\end{align*}
Hence,
\begin{equation}\label{TCase2}
    T^{1-2A(\log U)^{-\frac23}(\log\log U)^{-\frac13}-\delta_1} \leq T^{1-2^{\frac13}A(\delta_1\log T)^{-\frac23}(\log \log T)^{-\frac13}-\delta_1}.
\end{equation}
Also,
\begin{equation*}
    U \geq T^{\delta_1-\delta_1^2/2} \geq T^{\frac{\delta_1}{2}} \implies \log U \geq \frac12\delta_1\log T,
\end{equation*}
and
\begin{align*}
    \log \log U &\geq \log \left(\frac12\delta_1\log T\right) \geq \log\left(\frac12(\log T)^{\frac35}(\log\log T)^{-\frac15}\right)
    \\&\geq
    \log\left((\log T)^{\frac25}\right)=\frac25 \log\log T.
\end{align*}
Therefore,
\begin{align*}
    U&^{\frac32A(\log U)^{-\frac23}(\log\log U)^{-\frac13}+\frac{\delta_1}{2}} \leq U^{4A(\delta_1\log T)^{-\frac23}(\log\log T)^{-\frac13}+\frac{\delta_1}{2}}
    \\&=
    T^{\left((\log T)^{-\frac25}(\log\log T)^{-\frac15}+\delta_1-\frac{\delta_1^2}{2}\right)\left(4A(\delta_1\log T)^{-\frac23}(\log\log T)^{-\frac13}+\frac{\delta_1}{2}\right)}
    \\&\leq
    T^{o\left((\delta_1\log T)^{-\frac23}(\log\log T)^{-\frac13}\right)+4A\delta_1^{\frac13}(\log T)^{-\frac23}(\log\log T)^{-\frac13}+\frac{\delta_1}{2}(\log T)^{-\frac25}(\log\log T)^{-\frac15}+\frac{\delta_1^2}{2}-\frac{\delta_1^3}{4}}
    \\&\leq
    T^{o\left((\delta_1\log T)^{-\frac23}(\log\log T)^{-\frac13}\right)+\delta_1(\log T)^{-\frac25}(\log\log T)^{-\frac15}+\frac{\delta_1^2}{2}-\frac{\delta_1^3}{4}}.
\end{align*}
If $\delta_1 \geq (\log T)^{-\frac15}$, then $\delta_1(\log T)^{-\frac25}(\log\log T)^{-\frac15}<\frac{\delta_1^3}{4}$.\\
If $\delta_1 \leq (\log T)^{-\frac15}$, then $\delta_1(\log T)^{-\frac25}(\log\log T)^{-\frac15}=o\left((\delta_1\log T)^{-\frac23}(\log\log T)^{-\frac13}\right)$.\\
Thus,
\begin{equation}\label{UCase2}
    U^{\frac32A(\log U)^{-\frac23}(\log\log U)^{-\frac13}+\frac{\delta_1}{2}} \leq T^{\frac{\delta_1^2}{2}+o\left((\delta_1\log T)^{-\frac23}(\log\log T)^{-\frac13}\right)}.
\end{equation}
By \eqref{TCase2} and \eqref{UCase2}, we find that the right-hand side of \eqref{integral} is $\ll T^{1-A(\delta_1\log T)^{-\frac23}(\log\log T)^{-\frac13}-\delta_1+\delta_1^2/2}$.

\section{Proof of Theorem \ref{thm1}}\label{sec3}

Assuming the Riemann hypothesis and using the Taylor series, we have
    \begin{equation}\label{Sthm1}
        \mathcal{S}(\delta)=\sum_{0<\gamma \leq T}|\zeta(\rho+\delta)|^2=\sum_{0<\gamma \leq T}\left|\zeta'(\rho)\delta+\sum_{n=2}^{\infty}\frac{\zeta^{(n)}(\rho)}{n!}\delta^n\right|^2.
    \end{equation}
The Riemann hypothesis implies the Lindel\"of hypothesis, which states that uniformly in any fixed vertical strip with $t>2$ and for any $\epsilon>0$ we have (see \cite[Chapters 13 and 14]{Titchmarsh})
\begin{equation}\label{zetabound}
    \zeta(\sigma+it) \ll
    \begin{cases}
        t^{\epsilon} \quad \text{if} \quad \sigma \geq \frac12, \\
        t^{\frac12-\sigma+\epsilon} \quad \text{if} \quad \sigma \leq \frac12.
    \end{cases}
\end{equation}
By Cauchy's integral formula for derivatives, we have, for $n \geq 1$ and $0<\gamma \leq T$,
    \begin{equation}\label{CIFbound}
        \frac{\zeta^{(n)}(\rho)}{n!}=\frac{1}{2\pi i}\int\limits_{|s-\rho|=\frac{1}{\log \gamma}}\frac{\zeta(s)}{(s-\rho)^{n+1}}ds \ll T^{\epsilon}\log^nT 
    \end{equation}
uniformly in $n$.  For large \(T\) we may assume that \(|\delta|\log T \le \tfrac{1}{2}\).
  Hence, for $0<\gamma \leq T$,
    \begin{equation*}
\sum_{n=2}^{\infty}\frac{\zeta^{(n)}(\rho)}{n!}\delta^n \ll T^{\epsilon} \sum_{n=2}^{\infty} (|\delta|\log T)^n \ll |\delta|^2 T^{\epsilon} \log^2 T
    \end{equation*}
uniformly in $\delta$. So, for $0<\gamma \leq T$,
    \begin{equation*}
\sum_{n=2}^{\infty}\frac{\zeta^{(n)}(\rho)}{n!}\delta^n \ll |\delta|^2 T^{\epsilon}.
    \end{equation*}
In \cite{Gonek}, Gonek proved the formula
    \begin{equation*}
        \sum_{0<\gamma \leq T}\zeta'(\rho)\zeta'(1-\rho)=\frac{1}{12}\frac{T}{2\pi}\log^4T+O(T\log^3T).
    \end{equation*}
So, assuming RH,
    \begin{equation*}
        \sum_{0<\gamma \leq T}|\zeta'(\rho)|^2=\frac{1}{12}\frac{T}{2\pi}\log^4T+O(T\log^3T).
    \end{equation*}
By \eqref{CIFbound}, for $0 < \gamma \leq T$ we have $\zeta'(\rho) \ll T^{\epsilon}\log T$, so the formula $N(T):=\sum\limits_{0<\gamma \leq T}1 \sim \frac{T}{2\pi}\log T$ leads to
\begin{equation*}
    \sum_{0<\gamma\le T}|\zeta'(\rho)| \ll T^{1+\epsilon}\log^2 T.
\end{equation*}
Hence, we conclude
    \begin{equation*}
    \begin{split}
        \mathcal{S}(\delta)&=\sum_{0<\gamma \leq T}\left|\zeta'(\rho)\delta+O(|\delta|^2 T^{\epsilon})\right|^2
        \\&=
        \frac{|\delta|^2}{12}\frac{T}{2\pi}\log^4T+O(|\delta|^2T\log^3T)+O\left(|\delta|^4T^{2\epsilon}\sum_{0<\gamma \leq T}1\right)+O\left(|\delta|^3T^{\epsilon}\sum_{0<\gamma \leq T}|\zeta'(\rho)|\right)
        \\&=
        \frac{|\delta|^2}{12}\frac{T}{2\pi}\log^4T+O(|\delta|^2T\log^3T)+O(|\delta|^3T^{1+2\epsilon}\log^2 T).
    \end{split}
    \end{equation*}
Therefore, for all $\epsilon>0$ we have
    \begin{equation*}
        \mathcal{S}(\delta)=\frac{|\delta|^2}{12}\frac{T}{2\pi}\log^4T+O(|\delta|^2T\log^3T)+O(|\delta|^3T^{1+\epsilon})
    \end{equation*}
uniformly in $\delta$.

\section{Proof of Theorem \ref{thm2}}\label{sec4}

Let $L=(2a+1)^{-1}\log \frac{T}{2\pi}$ and $c=1+L^{-1}$. Let $U$ be such that $U \to \infty$ as $T \to \infty$ and $U \ll T^{1-\epsilon}$ for some $\epsilon>0$. Assume that $-\frac{a}{\log T} \leq \delta_1 \leq \frac12+\frac{a}{\log T}$, $\delta_1 \neq 0$.

Recall (\cite[formulae (16) and (18)]{Garunkstis}) the functional equation $\zeta(s)=\Delta(s)\zeta(1-s)$ and the formula
\begin{equation}\label{Delta}
    \frac{\Delta'}{\Delta}(s)=\frac{i\pi}{2}-\log\frac{s}{2\pi}+g(s)\qquad(t>2),
\end{equation}
where $g(s)$ is analytic, $g(s) \ll |s|^{-1}$, and $g'(s) \ll |s|^{-2}$.

In \cite[formula (28)]{Garunkstis}, using Gonek's \cite{Gonek} approach it was shown that if $-\frac{a}{\log T} \leq \delta_1 \leq \frac12+\frac{a}{\log T}$ then
\begin{equation}\label{S}
    \mathcal{S}(\delta)=2 \Re (J_1)-J_2+O\left(\sqrt{T}\log^4T\right),
\end{equation}
where
\begin{align*}
    &J_1=\frac{1}{2\pi i}\int\limits_{c-iU}^{c+iU}\frac{\zeta'}{\zeta}(s+\delta)\zeta(s)\zeta(s+2\delta_1)\left(\frac{T}{2\pi}\right)^s\frac{ds}{s}+O\left(\sqrt{T}\log^3T\right)+O\left(\frac{T\log^3T}{U}\right),\\
    &J_2=\int\limits_1^T\frac{\Delta'}{\Delta}\left(-L^{-1}+it-\delta\right)d\mathcal{H}(t)+O\left(\sqrt{T}\log^3T\right),
\end{align*}
with
\begin{equation*}
\overline{\mathcal{H}(\tau)}=\frac{1}{2\pi i}\int\limits_{c-iU_{\tau}}^{c+iU_{\tau}}\zeta(s)\zeta(s+2\delta_1)\left(\frac{\tau}{2\pi}\right)^s\frac{ds}{s}+O\left(\sqrt{T}\log^2T\right)+O\left(\frac{\tau\log^2\tau}{U_{\tau}}\right).
\end{equation*}

To calculate $J_1$, take $U=T^{1-\epsilon}$.\\
If $d<1$ is such that $d+\delta_1>\frac12$ and $d\neq 1-\delta_1,1-2\delta_1$, then assuming the Riemann hypothesis, we have
\begin{align}\label{J1integral}
J_1
=&
\frac{\zeta'}{\zeta}(1+\delta)\zeta(1+2\delta_1)\frac{T}{2\pi} \nonumber
\\& +
\mathbbm{1}_{(1-c,1-d)}(2\delta_1)\frac{\zeta'}{\zeta}(1-\overline{\delta})\frac{\zeta(1-2\delta_1)}{1-2\delta_1}\left(\frac{T}{2\pi}\right)^{1-2\delta_1}\nonumber
\\&-
\mathbbm{1}_{(1-c,1-d)}(\delta_1)\zeta(1+\overline{\delta})\frac{\zeta(1-\delta)}{1-\delta}\left(\frac{T}{2\pi}\right)^{1-\delta}
    \\&-
    \frac{1}{2\pi i}\left(\int\limits_{c+iU}^{d+iU}+\int\limits_{d+iU}^{d-iU}+\int\limits_{d-iU}^{c-iU}\right)\frac{\zeta'}{\zeta}(s+\delta)\zeta(s)\zeta(s+2\delta_1)\left(\frac{T}{2\pi}\right)^s\frac{ds}{s}\nonumber
\\&+
O\left(T^{\frac12+\epsilon}\right).\nonumber
\end{align}
We can choose
\begin{equation*}
d=
    \begin{cases}
        \frac12 &\text{ if } \delta_1>\epsilon, \delta_1\neq \frac14, \delta_1 \neq \frac12,\\
        \frac12+\frac{2a}{\log T} &\text{ if } \delta_1\leq \epsilon \text{ or } \delta_1=\frac14 \text{ or } \delta_1=\frac12.
    \end{cases}
\end{equation*}
For $-1 \leq \sigma \leq 2, ~ 2\leq t \leq T$, we have (\cite[Theorem 9.6 (A)]{Titchmarsh})
\begin{equation*}
    \frac{\zeta'}{\zeta}(s)=\sum_{|t-\gamma|\leq 1}\frac{1}{s-\rho}+O(\log T).
\end{equation*}
Hence, assuming the Riemann hypothesis, for $d \leq \sigma \leq 2, ~ 2 \leq t \leq T$, we have
\begin{align*}
    \frac{\zeta'}{\zeta}(s+\delta)&=\sum_{|t-\gamma|\leq 1}\frac{1}{\sigma+\delta_1-\frac{1}{2}+i(t+\delta_2-\gamma)}+O(\log T)\\
    &\ll
    \sum_{|t-\gamma|\leq 1} \log T +O(\log T) \ll \log^2 T.
\end{align*}
So, 
\begin{equation}\label{bound}
    \frac{\zeta'}{\zeta}(s+\delta) \ll \log^2 T
\end{equation}
uniformly in $d \leq \sigma \leq c, ~ 2\leq t \leq U$.

By \eqref{zetabound} and \eqref{bound}, $\frac{\zeta'}{\zeta}(s+\delta)\zeta(s)\zeta(s+2\delta_1)T^s \ll T^{1+\epsilon}$ on $d \leq \sigma \leq c, ~ t=U$. Also, $1/s \ll U^{-1}=T^{-1+\epsilon}$ on $d \leq \sigma \leq c, ~ t=U$. Hence, for all $\epsilon>0$,
\begin{equation}\label{intbound1}
    \int\limits_{c \pm iU}^{d \pm iU}\frac{\zeta'}{\zeta}(s+\delta)\zeta(s)\zeta(s+2\delta_1)\left(\frac{T}{2\pi}\right)^s\frac{ds}{s} \ll T^{\epsilon}.
\end{equation}

By \eqref{zetabound} and \eqref{bound}, $\frac{\zeta'}{\zeta}(s+\delta)\zeta(s)\zeta(s+2\delta_1)T^s \ll T^{\frac12+\epsilon}$ on $\sigma=d, ~ 2 \leq t \leq U$. Therefore, for all $\epsilon>0$,
\begin{equation}\label{intbound2}
    \int\limits_{d-iU}^{d+iU} \frac{\zeta'}{\zeta}(s+\delta)\zeta(s)\zeta(s+2\delta_1)\left(\frac{T}{2\pi}\right)^s\frac{ds}{s} \ll T^{\frac{1}{2}+\epsilon}.
\end{equation}

By \eqref{J1integral}, \eqref{intbound1} and \eqref{intbound2}, we get
\begin{align*}
    J_1=&\frac{\zeta'}{\zeta}(1+\delta)\zeta(1+2\delta_1)\frac{T}{2\pi}+\mathbbm{1}_{(1-c,1-d)}(2\delta_1)\frac{\zeta'}{\zeta}(1-\overline{\delta})\frac{\zeta(1-2\delta_1)}{1-2\delta_1}\left(\frac{T}{2\pi}\right)^{1-2\delta_1}
    \\&-\mathbbm{1}_{(1-c,1-d)}(\delta_1)\zeta(1+\overline{\delta})\frac{\zeta(1-\delta)}{1-\delta}\left(\frac{T}{2\pi}\right)^{1-\delta}+O\left(T^{\frac{1}{2}+\epsilon}\right).\nonumber
\end{align*}

We have $1-c=-\frac{2a+1}{\log\frac{T}{2\pi}}<-\frac{a}{\log T}<\delta_1$ and $1-c=-\frac{2a+1}{\log\frac{T}{2\pi}}<-\frac{2a}{\log T}<2\delta_1$.

If $\delta_1>1-d$, then $T^{1-\delta}\leq T^d\ll\sqrt{T}$ and $\zeta(1+\overline{\delta})\frac{\zeta(1-\delta)}{1-\delta} \ll 1$.

Thus,
\begin{align*}
    J_1=&\frac{\zeta'}{\zeta}(1+\delta)\zeta(1+2\delta_1)\frac{T}{2\pi}+\mathbbm{1}_{\left[-\frac{a}{\log T},\frac{1-d}{2}\right)}(\delta_1)\frac{\zeta'}{\zeta}(1-\overline{\delta})\frac{\zeta(1-2\delta_1)}{1-2\delta_1}\left(\frac{T}{2\pi}\right)^{1-2\delta_1}
    \\&-\zeta(1+\overline{\delta})\frac{\zeta(1-\delta)}{1-\delta}\left(\frac{T}{2\pi}\right)^{1-\delta}+O\left(T^{\frac{1}{2}+\epsilon}\right).\nonumber
\end{align*}

We have $\frac{1-d}{2}=\frac14 \text{ or }\frac14-\frac{a}{\log T}$.\\
If $\frac14-\frac{a}{\log T} \leq \delta_1 \leq \frac14$, then $\frac{\zeta'}{\zeta}(1-\overline{\delta})\frac{\zeta(1-2\delta_1)}{1-2\delta_1}\left(\frac{T}{2\pi}\right)^{1-2\delta_1}\ll \sqrt{T}$. Hence,

\begin{align}\label{J1}
    J_1=&\frac{\zeta'}{\zeta}(1+\delta)\zeta(1+2\delta_1)\frac{T}{2\pi}+\mathbbm{1}_{\left[-\frac{a}{\log T},\frac14\right)}(\delta_1)\frac{\zeta'}{\zeta}(1-\overline{\delta})\frac{\zeta(1-2\delta_1)}{1-2\delta_1}\left(\frac{T}{2\pi}\right)^{1-2\delta_1}
    \\&-\zeta(1+\overline{\delta})\frac{\zeta(1-\delta)}{1-\delta}\left(\frac{T}{2\pi}\right)^{1-\delta}+O\left(T^{\frac{1}{2}+\epsilon}\right).\nonumber
\end{align}

Now we will calculate $J_2$.
Let $U_{\tau}=\tau^{1-\epsilon}$. If $d<1$ is such that $d+\delta_1>\frac12$ and $d\neq 1-2\delta_1$, then assuming the Riemann hypothesis, we obtain
\begin{align*}
&\overline{\mathcal{H}(\tau)}=\zeta(1+2\delta_1)\frac{\tau}{2\pi}+\mathbbm{1}_{(1-c,1-d)}(2\delta_1)\frac{\zeta(1-2\delta_1)}{1-2\delta_1}\left(\frac{\tau}{2\pi}\right)^{1-2\delta_1}
\\&-\frac{1}{2\pi i}\left(\int\limits_{c+iU_{\tau}}^{d+iU_{\tau}}+\int\limits_{d+iU_{\tau}}^{d-iU_{\tau}}+\int\limits_{d-iU_{\tau}}^{c-iU_{\tau}}\right)\zeta(s)\zeta(s+2\delta_1)\left(\frac{\tau}{2\pi}\right)^s\frac{ds}{s}+O\left(T^{\frac12+\epsilon}\right).
\end{align*}

We can choose
\begin{equation*}
d=
    \begin{cases}
        \frac12 &\text{ if } \delta_1>0, \delta_1 \neq \frac14,\\
        \frac12+\frac{2a}{\log T} &\text{ if } \delta_1\leq 0 \text{ or } \delta_1=\frac14.
    \end{cases}
\end{equation*}

By \eqref{zetabound} for $1\leq \tau \leq T$, we have
\begin{equation*}
    \int\limits_{c \pm iU_{\tau}}^{d \pm iU_{\tau}}\zeta(s)\zeta(s+2\delta_1)\left(\frac{\tau}{2\pi}\right)^s\frac{ds}{s} \ll T^{\epsilon},
\end{equation*}
and
\begin{equation*}
    \int\limits_{d+2i}^{d+iU_{\tau}}\zeta(s)\zeta(s+2\delta_1)\left(\frac{\tau}{2\pi}\right)^s\frac{ds}{s} \ll T^{\frac{1}{2}+\epsilon}.
\end{equation*}
Also,
\begin{equation*}
    \frac{\zeta(1-2\delta_1)}{1-2\delta_1}\left(\frac{\tau}{2\pi}\right)^{1-2\delta_1}\ll \sqrt{T}
\end{equation*}
uniformly in $\frac14-\frac{a}{\log T} \leq \delta_1 \leq \frac14$, $1 \leq \tau \leq T$.
Hence, for $1\leq \tau \leq T$,
\begin{equation*}
    \mathcal{H}(\tau)=\zeta(1+2\delta_1)\frac{\tau}{2\pi}+\mathbbm{1}_{\left[-\frac{a}{\log T},\frac14\right)}(\delta_1)\frac{\zeta(1-2\delta_1)}{1-2\delta_1}\left(\frac{\tau}{2\pi}\right)^{1-2\delta_1}+O\left(T^{\frac{1}{2}+\epsilon}\right).
\end{equation*}
So, for $1\leq t \leq T$,
\begin{equation*}
    \mathcal{H}(t)=\zeta(1+2\delta_1)\frac{t}{2\pi}+\mathbbm{1}_{\left[-\frac{a}{\log T},\frac14\right)}(\delta_1)\frac{\zeta(1-2\delta_1)}{1-2\delta_1}\left(\frac{t}{2\pi}\right)^{1-2\delta_1}+h(t),
\end{equation*}
where $h(t)=O\left(T^{\frac{1}{2}+\epsilon}\right)$.

In the same way as in \cite{Garunkstis}, it can be shown that
\begin{align*}
    \int\limits_1^T\frac{\Delta'}{\Delta}&\left(-L^{-1}+it-\delta\right)d(\mathcal{H}(t)-h(t))=-\zeta(1+2\delta_1)\frac{T}{2\pi}\left(\log\frac{T}{2\pi}-1\right)
    \\&-
    \mathbbm{1}_{\left[-\frac{a}{\log T},\frac14\right)}(\delta_1)\frac{\zeta(1-2\delta_1)}{1-2\delta_1}\left(\frac{T}{2\pi}\right)^{1-2\delta_1}\left(\log\frac{T}{2\pi}-\frac{1}{1-2\delta_1}\right)+O\left(\log^2T\right).
\end{align*}

By integration by parts and using \eqref{Delta}, we get
\begin{align*}
    \int\limits_1^T&\frac{\Delta'}{\Delta}\left(-L^{-1}+it-\delta\right)dh(t)
    \\=&
    \left[\left(\frac{i\pi}{2}-\log\frac{-L^{-1}+it-\delta}{2\pi}+g(-L^{-1}+it-\delta)\right)h(t)\right]_1^T
    \\&-
    \int\limits_1^Th(t)d\left(\frac{i\pi}{2}-\log\frac{-L^{-1}+it-\delta}{2\pi}+g(-L^{-1}+it-\delta)\right)
    \\\ll&
    T^{\frac12+\epsilon}.
\end{align*}

Hence,
\begin{align}\label{J2}
    J_2=-&\zeta(1+2\delta_1)\frac{T}{2\pi}\left(\log\frac{T}{2\pi}-1\right)\nonumber
    \\&-
    \mathbbm{1}_{\left[-\frac{a}{\log T},\frac14\right)}(\delta_1)\frac{\zeta(1-2\delta_1)}{1-2\delta_1}\left(\frac{T}{2\pi}\right)^{1-2\delta_1}\left(\log\frac{T}{2\pi}-\frac{1}{1-2\delta_1}\right)
    \\&+
    O\left(T^{\frac{1}{2}+\epsilon}\right).\nonumber
\end{align}

By \eqref{S}, \eqref{J1}, and \eqref{J2} we obtain
\begin{align*}
    \mathcal{S}(\delta)=&\zeta(1+2\delta_1)\frac{T}{2\pi}\left(\log\frac{T}{2\pi}-1\right)+\mathbbm{1}_{\left[-\frac{a}{\log T},\frac14\right)}(\delta_1)\frac{\zeta(1-2\delta_1)}{1-2\delta_1}\left(\frac{T}{2\pi}\right)^{1-2\delta_1}\left(\log\frac{T}{2\pi}-\frac{1}{1-2\delta_1}\right)\nonumber
    \\&+
    2\Re\left(\frac{\zeta'}{\zeta}(1+\delta)\zeta(1+2\delta_1)\frac{T}{2\pi}\right)+2\mathbbm{1}_{\left[-\frac{a}{\log T},\frac14\right)}(\delta_1)\Re\left(\frac{\zeta'}{\zeta}(1-\overline{\delta})\frac{\zeta(1-2\delta_1)}{1-2\delta_1}\left(\frac{T}{2\pi}\right)^{1-2\delta_1}\right)
    \\&-2\Re\left(\zeta(1+\overline{\delta})\frac{\zeta(1-\delta)}{1-\delta}\left(\frac{T}{2\pi}\right)^{1-\delta}\right)
    +
    O\left(T^{\frac{1}{2}+\epsilon}\right).\nonumber
\end{align*}

This completes the proof of Theorem $\ref{thm2}$.

\section{Asymptotics. Proof of Corollary \ref{cor}}\label{sec5}

Write \eqref{Sthm2} as
\begin{equation*}
    \mathcal{S}(\delta)=\frac{T}{2\pi}\Re g(\delta)+O\left(T^{\frac{1}{2}+\epsilon}\right),
\end{equation*}
where
\begin{align*}
g(\delta)=&\zeta(1+2\delta_1)\log\frac{T}{2\pi e}+\mathbbm{1}_{\left[-\frac{a}{\log T},\frac14\right)}(\delta_1)\frac{\zeta(1-2\delta_1)}{1-2\delta_1}\left(\frac{T}{2\pi}\right)^{-2\delta_1}\left(\log\frac{T}{2\pi}-\frac{1}{1-2\delta_1}\right)
\\&+
2\frac{\zeta'}{\zeta}(1+\delta)\zeta(1+2\delta_1)+2\mathbbm{1}_{\left[-\frac{a}{\log T},\frac14\right)}(\delta_1)\frac{\zeta'}{\zeta}(1-\overline{\delta})\frac{\zeta(1-2\delta_1)}{1-2\delta_1}\left(\frac{T}{2\pi}\right)^{-2\delta_1}
\\&-
2\zeta(1+\overline{\delta})\frac{\zeta(1-\delta)}{1-\delta}\left(\frac{T}{2\pi}\right)^{-\delta}.
\end{align*}

In the proof of Theorem 6 in \cite{Garunkstis} it was shown that if $|\delta|\log T \leq 1$ then
\begin{equation*}
    \Re g(\delta)=\frac{|\delta|^2}{12}\log^4T+O\left(|\delta|^2\log^3T\right)+O\left(|\delta|^3\log^5T\right).
\end{equation*}

From this we see that there is $C>0$ such that $\Re g(\delta) \geq \frac{|\delta|^2}{12}\log^4T-C|\delta|^2\log^3T-C|\delta|^3\log^5T$ in $|\delta|\log T \leq 1$. Hence, there is $\epsilon_0>0$ such that $\Re g(\delta) \gg |\delta|^2\log^4T$ uniformly in $|\delta|\log T \leq \epsilon_0$. So,
\begin{equation*}
    \frac{T^{\frac12+\epsilon}}{T\Re g(\delta)} \ll \frac{1}{T^{\frac12-\epsilon}|\delta|^2} \to 0 \quad \text{if} \quad |\delta|\log T \leq \epsilon_0 ~ \text{and} ~ |\delta|T^{\frac14-\epsilon} \to \infty.
\end{equation*}

Hence, if $\delta_1 \neq 0$, $|\delta|\log T \leq \epsilon_0$ and $|\delta| \geq T^{-\frac14+\epsilon}$, then
\begin{equation*}
S(\delta)\sim\frac{T}{2\pi}\Re g(\delta).
\end{equation*}

\bigskip

\bibliographystyle{plain}
\bibliography{Discrete}

\end{document}